\documentclass[12pt]{article}
\usepackage{amscd,amsfonts,amsmath,amssymb,amstext,amsthm}

\title{On closures of cycle spaces of flag domains}
\author{Jaehyun Hong and Alan Huckleberry}
\date{\today}
\def\dim{\operatorname{dim}}
\def\exp{\operatorname{exp}}
\def\cl{\operatorname{cl}}
\def\bd{\operatorname{bd}}
\theoremstyle{plain}
\newtheorem{theorem} {Theorem} [section]
\newtheorem{lemma} [theorem]{Lemma}
\newtheorem{proposition}[theorem]{Proposition}

\theoremstyle{definition}

\begin{document}
\maketitle
\abstract{\begin{quote} {\small Open orbits
$D$ of noncompact real forms $G_0$ acting on flag manifolds
$Z=G/Q$ of their semisimple complexifications $G$ are considered.
Given $D$ and a maximal compact subgroup $K_0$ of $G_0$, there is
a unique complex $K_0$--orbit in $D$ which is regarded as a point
$C_0\in {\mathcal C}_q(D)$ in the space of $q$-dimensional cycles
in $D$.  The group theoretical cycle space $\mathcal M_D$ is
defined to be the connected component containing $C_0$ of the
intersection of the $G$--orbit $G(C_0)$ with ${\mathcal C}_q(D)$.
The main result of the present article is that ${\mathcal M}_D$ is
closed in ${\mathcal C}_q(D)$.  This follows from an analysis of
the closure of the universal domain $\mathcal U$ in any
$G$-equivariant compactification of the affine symmetric space
$G/K$, where $K$ is the complexification of $K_0$ in $G$.}
\end {quote}}

\section {Background and notation}
\label {introduction}

Throughout this article $G_0$ denotes a noncompact simple Lie
group of adjoint type.  Generalizations of our results to the
semisimple case require only formal adjustments and will not be
discussed. Having fixed a maximal compact subgroup $K_0$ of $G_0$,
we will make use
of Iwasawa decompositions $G_0=K_0A_0N_0$.\\

We regard $G_0$ as a closed subgroup of its universal
complexification $G$.   Therefore the complexification $K$ of
$K_0$ is a closed complex subgroup of $G$ and we consider the affine
symmetric space $\Omega =G/K$.  Certain $G$-equivariant
projective algebraic compactifications $X$ of $\Omega $ play an
important role in our work.
These arise as follows.\\

Let $Z=G/Q$ be a $G$-flag manifold, i.e., $Q$ is a complex
parabolic subgroup of $G$.  Starting with the basic work (\cite
{W1}) there has been substantial interest in complex geometric
objects related to the $G_0$-action on $Z$ (see \cite {FHW} for a
systematic presentation). In particular there are only finitely
many $G_0$--orbits in $Z$ and therefore there are open orbits $D$
which merit study from the
complex geometric viewpoint.\\

Each such $D$ contains a unique $K_0$--orbit $C_0$ which is a
\emph{complex} submanifold of $Z$.  We let $q:=\dim_\mathbb CC_0$
be the dimension of this \emph{base cycle} and regard it as a
point $C_0\in {\mathcal C}_q(D)$ in the full cycle space of $D$.
Without further notation we replace ${\mathcal C}_q(D)$ by its
connected component containing $C_0$ in the irreducible component
which contains $C_0$ in the full cycle space ${\mathcal C}_q(Z)$.\\

In (\cite {WeW}) a group theoretical cycle space ${\mathcal M}_D$
was introduced.  For this consider the $G$-orbit ${\mathcal
M}_Z:=G(C_0)$ of the base cycle in ${\mathcal C}_q(Z)$. Since the
induced $G$--action on ${\mathcal C}_q(Z)$ is algebraic,
${\mathcal M}_Z$ is Zariski open in its closure and its
intersection ${\mathcal M}_Z\cap {\mathcal C}_q(D)$ with the
semialgebraic open set ${\mathcal C}_q(D)$ consists of at most
finitely many components. The cycle space ${\mathcal M}_D$ is
defined to be the connected component
of this intersection which contains $C_0$.

\section {The main result}
Our goal here is to present a proof of the following result. It
will be reliant on the more technical results of the following
sections.
\begin {theorem}\label{M_D is closed}
The group theoretical cycle space ${\mathcal M}_D$ is closed in
the full cycle space ${\mathcal C}_q(D)$.
\end {theorem}
Let us attempt to put this in perspective.  First
of all ${\mathcal M}_D$ is a locally closed complex submanifold
of ${\mathcal C}_q(D)$.  The representation of the $G$--isotropy
group on the tangent space of ${\mathcal C}_q(D)$ at the
base cycle $C_0$ has been calculated in detail (see Part IV
in \cite{FHW}).  In particular, even for a fixed $G_0$, there
is a great variety of representations depending on $D$ and
the flag manifold $Z$, and, for example, the codimension
of ${\mathcal M}_D$ in ${\mathcal C}_q(D)$ can vary wildly.\\

In the case where ${\mathcal M}_D$ is open in ${\mathcal C}_q(D)$
Theorem \ref{M_D is closed} states that 
${\mathcal M}_D={\mathcal C}_q(D)$.  This is particularly
useful in situations where the cycles have meaning in 
complex geometry, e.g., in the case of period domains
such as the moduli space of marked K3-surfaces: 
Any two cycles differ only by a transformation in the
complex group $G$ and no degeneration is possible.

In the future we hope that the full cycle space ${\mathcal C}_q(D)$
can be explicitly computed and that it will be of use
in representation theory.  We already know that in many
cases it is a Stein space and it is very likely that it
is Kobayashi hyperbolic.  Since ${\mathcal M}_D$
is closed, it is quite possible that recently developed
$G_0$--invariant theory (see \cite{HSch, HSt})
can be applied to show that ${\mathcal C}_q(D)$ is
a Luna-slice type bundle over ${\mathcal M}_D$.  Since
${\mathcal M}_D$ has already been described with great
precision (see below) and there are good Ans\"atze for
describing the fiber, this could very well lead to
the desired precise description of ${\mathcal C}_q(D)$.\\

Now let us recall the description of ${\mathcal M}_D$.
For a certain well--understood special class of 
domains which are said to be of 
\emph{Hermitian holomorphic type}, where in
particular $G_0$ is of Hermitian type, $\mathcal M_D$ is just
the associated bounded symmetric domain ${\mathcal B}$.  In this
case the stabilizer of $C_0$ in $G$ is a parabolic group
$P$ so that ${\mathcal M}_Z=G/P$ is the compact dual of  
${\mathcal B}$ (see e.g. \cite{FHW}).  Since 
${\mathcal M}_Z$ compact, 
it is a direct consequence of
the definitions that ${\mathcal M}_D$ is closed in 
${\mathcal C}_q(D)$.  Thus we may assume that $D$ is not of Hermitian holomorphic type, and the following result is applicable (\cite{{HW},{FH}}, see also \cite {FHW}).
\begin {theorem}\label{classification}
If $D$ is not of Hermitian holomorphic type, then ${\mathcal M}_D$
is naturally biholomorphic to the universal domain ${\mathcal U}$
contained in the affine symmetric space $\Omega =G/K$.
\end {theorem}
Before going into the details of the definition of $\mathcal U$
(first introduced in a representation theoretical context by
Akhiezer and Gindikin), we emphasize that $\mathcal U$ is defined
independent of which flag manifold $Z$ and domain $D$ is under
consideration.  For this reason and since $\mathcal U$ occurs in a
number of contexts, several of which are important in this
article, we refer to it as being
\emph{universal}.\\

The domain $\mathcal U$ is defined as follows.  If ${\mathfrak
g}_0={\mathfrak k}_0\oplus {\mathfrak a}_0 \oplus {\mathfrak n}_0$
is an Iwasawa decomposition at the Lie algebra level, then one
defines the polytope $\omega _0\subset {\mathfrak a}_0$ as
follows:
\begin {equation*}
\omega _0=\underset{\alpha}{\bigcap} \{\xi \in {\mathfrak
a}_0:\vert \alpha (\xi) \vert < \frac{\pi}{2}\}\,,
\end {equation*}
where $\alpha $ runs over the restricted roots. Then
\begin {equation*}
\mathcal U:=G_0\exp (i\omega _0)(x_0)\, ,
\end {equation*}
where $x_0$ is a base point with isotropy group $K$.\\

Since we have eliminated the Hermitian holomorphic case from
discussion, the $G$-isotropy group $\tilde K$ at $C_0$ is just a
finite extension of the connected group $K$, and the cycle
space $\mathcal M_D$ lifts biholomorphically to $\mathcal U$ in
$\Omega =G/K$ (\cite {FHW}).  In our discussion of closures
the finite cover $G/K\to G/\tilde K$ plays no
role. Hence, for notational 
convenience we simply assume that 
${\mathcal M}_D\subset \Omega $.
Since it is quite difficult to know anything specific about the
closure of ${\mathcal M}_Z$, we consider all possible situations
\begin {equation*}
\mathcal M_D=\mathcal U\subset \Omega \subset X=\cl(\Omega)\, ,
\end {equation*}
where $X$ is an arbitrary projective algebraic $G$-equivariant
compactification of $\Omega =G/K$.\\

One of the main methods used in proving Theorem
\ref{classification} is that of Schubert incidence geometry.
We make strong use of this in the present article, and therefore
we now sketch the basics (for details see \cite{FHW} or the original
papers \cite{{HW},{FH}}).
Given a Borel subgroup $B\subset G$, a $B$--Schubert variety $S$
in $Z$ is the closure $S=\mathcal O\dot \cup Y$ of a $B$-orbit
$\mathcal O$ in $Z$.  Here $q$-codimensional Schubert varieties of
Iwasawa-Borel subgroups, i.e., those which contain a component
$A_0N_0$ of an Iwasawa-decomposition $G_0=K_0A_0N_0$, are
important.\\

If $B$ is an Iwasawa-Borel subgroup and $S$ is an associated
$q$-codimensional Schubert variety with $S\cap C_0\not
=\emptyset$, then $S\cap C\not =\emptyset $ for every $C\in
\mathcal M_D$. Furthermore, the complement $Y$ of the open orbit
$\mathcal O$ in $S$ is contained in the complement of $D$.  In
particular, the incidence variety
\begin {equation*}
I_Y:=\{C\in \mathcal C_q(Z):C\cap Y\not =\emptyset \}
\end {equation*}
is contained in the complement of the cycle space $\mathcal M_D$.\\

The intersection $H_{Y, \Omega}:=I_Y\cap \Omega $ is a
$B$-invariant complex algebraic hypersurface which is contained in
the complement of $\mathcal M_D=\mathcal U$ in $\Omega $.  A final
step in the proof of Theorem \ref{classification} can be
formulated as follows. For this, in the Hermitian case the
Schubert variety $S$ must be chosen appropriately, but otherwise
it only required to have the properties described above.
\begin {theorem}\label {incidence hypersurfaces}
For every point $p$ in the boundary $\bd_\Omega(\mathcal
M_D)=\bd_\Omega(\mathcal U)$ there exists $k\in K_0$ with $p\in
k(H_{Y, \Omega})$.
\end {theorem}
We now state our main technical result which is proved in
$\S\ref{Iwasawa-envelope}$.
\begin {theorem} \label{no new interior}
Let $X$ be an arbitrary $G$--equivariant compactification of the
affine symmetric space $\Omega =G/K$.  Then the interior of the closure
$\cl _X({\mathcal U})$ of the universal domain ${\mathcal U}$ in $X$ is
${\mathcal U}$ itself and
\begin {equation*}
\bd _X({\mathcal U})=\cl _X(\bd _\Omega ({\mathcal U}))\,.
\end {equation*}
\end {theorem}
The essential point of this theorem is that the interior
of $\cl _X({\mathcal U})$ is ${\mathcal U}$.  The second statement,
$\bd _X({\mathcal U})=\cl _X(\bd _\Omega ({\mathcal U}))$, is a direct
consequence of this fact.  To see this, note that 
$\cl _X({\mathcal U})\setminus \cl _X(\bd _\Omega ({\mathcal U}))$ is
open in $\cl _X({\mathcal U})$ and contains ${\mathcal U}$.
Therefore, by the first statement
in the theorem
this open set is exactly ${\mathcal U}$ and we have have the decomposition 
\begin  {equation*}
\cl _X({\mathcal U})={\mathcal U}\,\dot \cup \, 
\cl _X(\bd _\Omega ({\mathcal U}))
\end  {equation*}
which is equivalent to the desired result.\\

Using Theorem~\ref{no new interior}, we can now give the\\
 
\noindent 
{\it Proof of Theorem \ref{M_D is closed}.} 
Applying Theorem~\ref{no new interior}, if
$H_X$ denotes the closure in $X$ of the hypersurface $H_{Y, \Omega}$ of
Theorem \ref{incidence hypersurfaces}, then, since $K_0$ is
compact, it follows that every $p\in \bd _X(\mathcal U)$ is
contained in some translate $k(H_X)$.

Since $Y$ is closed, every cycle in $k(H_X)$ also has nonempty
intersection with $k(Y)$.  Thus if $p\in \bd_X(\mathcal U)$ is
regarded as a cycle $C$, then $C\not \subset D$.  Therefore
$\mathcal M_D=\mathcal U$ is
closed in $\mathcal C_q(D)$.\qed \\

From this proof one sees that the main new ingredients for this
result are to be found in Theorem \ref{no new interior}.  The
proof of this result is in turn heavily reliant on particular
properties of special $G$-equivariant compactifictions $X$ of
$\Omega$.  In the Hermitian case we choose $X=X_+\times X_-$ 
to be the product of the two compact Hermitian symmetric spaces (see
$\S$\ref{hermitian case}) and in the nonhermitian case we make
strong use of the DeConcini-Procesi compactification (see
$\S$\ref{extending sections} and $\S$\ref{Iwasawa-envelope}). 

The desired result for an arbitrary equivariant compactification
of $\Omega $ follows from the fact that any two such 
compactifications are equivariantly birationally equivalent.
In the case where ${\mathcal M}_D$ is actually
contained in a finite (algebraic) quotient 
$\tilde \Omega = G/\tilde K$ of $\Omega $, the quotient
map extends to an equivariant rational map of the special
compactification of $\Omega $ under consideration to
the closure $\tilde X$ of ${\mathcal M}_Z$ in the cycle space
${\mathcal C}_q(Z)$.  The arguments that show that
the birational maps which arise from the various compactifications
of $\Omega $ play no role in the discussion show that
such generically finite rational maps also play no role.  Thus, as stated
above, we simply assume that ${\mathcal M}_Z=\Omega =G/K$ 
from the beginning.

\section {The Hermitian case} \label {hermitian case}
In this section it is assumed that $G_0$ is of Hermitian type. 
In this case the parabolic subgroups of $G$ containing $K_0$ are $P_+=KS_+$ and $P_-=KS_-$,where $S_+$ and $S_-$ are unipotent part of $P_+$ and $P_-$. They
correspond to Hermitian symmetric spaces $X_+=G/P_+$ and
$X_-=G/P_-$ with the base point $x_+$ and $x_-$. Consider the
diagonal action of $G$ on $X_+ \times X_-$. Then the isotropy
group at $x=(x_+, x_-)$ is $K$ and the  affine symmetric space
$\Omega=G(x_+, x_-)=G/K$ is open dense in $X_+ \times X_-$. Write
$E=X_+ \times X_- \backslash \Omega$. It is known that the
universal domain $\mathcal U$ is equal to $\mathcal B_+ \times
\mathcal B_-=G_0x_+ \times G_0 x_-$. The following Lemma is proved
in the proof of Theorem 3.8 of \cite{WZ}.
\begin{lemma} \label{boundary(B)}
$ \bd_{X_+}(\mathcal B_+) \times \mathcal B_-$ and  $ \mathcal B_+
\times \bd_{X_-} (\mathcal B_-)$ are contained in $G/K=G(x_+,
x_-)$.
\end{lemma}
\begin{proposition} \label{Hermitian special}
For any $G$-orbit $\mathcal O$ in $E$, $\cl_{X_+ \times
X_-}(\mathcal B_+ \times \mathcal B_-) \cap \mathcal O$ has no
interior point in $\mathcal O$.
\end{proposition}
\begin{proof}
The boundary $\bd_{X_+ \times X_-}(\mathcal B_+ \times \mathcal
B_-)$ is the union  $ (\bd_{X_+}(\mathcal B_+) \times \mathcal
B_-) \cup ( \mathcal B_+ \times \bd_{X_-}( \mathcal B_-)) \cup
(\bd_{X_+}(\mathcal B_+) \times \bd_{X_-}(\mathcal B_-))$. By
Lemma \ref{boundary(B)} the first two subsets  $
\bd_{X_+}(\mathcal B_+) \times \mathcal B_-$ and $ \mathcal B_+
\times \bd_{X_-} (\mathcal B_-)$ are contained in $G/K=G(x_+,
x_-)$. So $\cl_{X_+ \times X_-}(\mathcal B_+ \times \mathcal B_-)
\cap \mathcal O$ is contained in $\bd_{X_+}(\mathcal B_+) \times
\bd_{X_-}(\mathcal B_-)$.

   If $\cl_{X_+ \times X_-}(\mathcal B_+ \times \mathcal B_-) \cap
\mathcal O$ has an interior point in $\mathcal O$, then the image
$\pi(\cl_{X_+ \times X_-}(\mathcal B_+ \times \mathcal B_-) \cap
\mathcal O)$ under the projection $\pi: \mathcal O \rightarrow
X_+$ would have an interior point in $X_+$ because $\pi$ is
$G$-equivariant and surjective. But $\pi(\cl_{X_+ \times
X_-}(\mathcal B_+ \times \mathcal B_-) \cap \mathcal O)$ is
contained in $\pi(\bd_{X_+}(\mathcal B_+ )\times
\bd_{X_-}(\mathcal B_-))=\bd_{X_+}(\mathcal B_+ )$ which has no
interior point in $X_+$.
\end{proof}
Let us now turn to the\\ 

\noindent
{\it Proof of Theorem~\ref{no new interior} in the case where $G_0$
is of Hermitian type.} Let $E=X \backslash \Omega$. We will show that 
$\cl_X(\mathcal U) \cap E$ has no interior point in $E$. Let $X_0=X_+
\times X_-$ be the $G$-equivariant compactification of $\Omega $
considered above and put $E_0=X_0 \backslash \Omega$. Since $X$ is
$G$-equivariantly birationally equivalent to $X_0$, we have the
following diagram:

\begin {gather*}
\begin {matrix}
 & \ \tilde X & \\
\  \ \ \swarrow & &\searrow \\
X_0 & & \ \ \ \ X
\end {matrix}
\end {gather*}
Let $\pi :\tilde X\to X_0$ and $p:\tilde X\to X$ denote the
respective proper modifications.\\

Assume that $\bd_X(\mathcal U) \cap E$ has an interior point in
$E$. Then $\bd_X(\mathcal U) \cap \mathcal O$ has an interior point
in $\mathcal O$ for some $G$-orbit $\mathcal O$ in $E$ of
codimension $1$ in $X$. Therefore the restriction $p: \hat{\mathcal
O}:=p^{-1}(\mathcal O) \rightarrow \mathcal O$ is biholomorphic,
because the indeterminant locus of $p$ has codimension $\geq 2$.
The other projection 
$\pi\vert \hat{\mathcal O}: \hat{\mathcal O} \rightarrow
\pi(\hat{\mathcal O})$ may have positive dimensional fibers,
but still $\pi(\hat{\mathcal
O})$ is a $G$-orbit $\mathcal O_0$ in $X_0$. Since 
$\pi\vert \hat{\mathcal O} $ is an
open map, $\pi(\bd_X(\mathcal U) \cap \mathcal O) \subset
\bd_{X_0}(\mathcal U) \cap \mathcal O_0 $ has an interior point in
$\mathcal O_0$, contrary to Proposition \ref{Hermitian special}
\hfill \qed 

\section {Non-Hermitian case}
For the remainder of this paper we assume that $G_0$ is
not of Hermitian type. Our work is devoted to proving 
Theorem~\ref{no new interior} in that case. 
Here we let $X:=X^W$ be the DeConcini-Procesi compactification of
$\Omega =G/K$ (\cite {DeCP}).  Orbits in the boundary
$E:=X^W\setminus \Omega $ are denoted by $O_I$, where $I$ is a
subset of $\{1,\ldots r\}$, and $S_I:=\cl (O_I)$. Recall
that for every such $I$ the compactification $X^W$ is realized in
$\mathbb P(V_I)\times P(V_J)$  where $J=\{1, \cdots, r\}
\backslash I$, and that the projection on the first factor defines
a $G$-equivariant morphism $\pi _I:X\to \mathbb P(V_I)$ to the
projective space of
the irreducible representation space $V_I$.

The restriction $\pi_I \vert S_I:S_I\to C_I=G/P_I$ is a fiber bundle
whose fiber is the DeConcini-Procesi compactification of an affine
symmetric space of a root theoretically distinguished Levi-factor
of $P_I$.

Since we have supposed that $G_0$ is not of Hermitian type, the image
$X_I:=Im(\pi _I)$ is another $G$-equivariant compactification of
the affine symmetric space $G/K$. In this case we denote by
$\Omega _I$ the open $G$-orbit in $X_I$. \\

\subsection  {Extending sections}\label {extending sections}

Let $H_I$ be the restriction to $X_I$ of the
hyperplane bundle $H$ of $\mathbb P(V_I)$.  Since the
vector space of sections of $H$ is an irreducible
$G$-representation space, it follows that the restriction map
defines isomorphisms
\begin {equation*}
\Gamma (\mathbb P(V_I),H)\cong \Gamma (X_I,H_I)
\cong \Gamma (C_I,H\vert C_I)\,.
\end {equation*}
If $L_I:=\pi _I^*H_I$, then, using the
isomorphism $\Gamma (S_I,L_I\vert S_I)\cong \Gamma
(C_I,H_I\vert C_I)$, we see that the restriction map
\begin {equation*}
R_I:\Gamma (X,L_I)\to \Gamma (S_I,L_I\vert S_I)
\end {equation*}
is surjective. We note that $R_I$ is $G$-equivariant.  Therefore
we may choose a $G$-invariant irreducible representation subspace
of $\Gamma (X,L_I)$ which is mapped isomorphically onto
$\Gamma (S_I,L_I\vert S_I)$.  In particular, if $B$ is a
Borel subgroup of $G$ and $s_0$ is a $B$-eigenvector in $\Gamma
(S_I,L_I\vert S_I)$, then there is a $B$-eigenvector $t_0\in
\Gamma (X,L_I)$ with $R_I(t_0)=s_0$.
\begin {proposition}\label {restriction} 
If $s \not=0\in \Gamma (S_I,L_I\vert S_I)$ is the restriction 
$s=R_I(t)$, then the intersection of the support $\vert t\vert $
with $\Omega $ is not empty.
\end {proposition}
\begin {proof}
If not, then $\vert t\vert $ is the union of certain irreducible
components of $X\setminus \Omega $.  But $t\vert S_I=\pi _I^*(t_I)$ and
$t_I$ is the restriction to $C_I$ of a unique section $\tilde
t_I\in \Gamma (X_I,H_I)$.  However, $\pi _I^*(\tilde
t_I)=t$, and since $\pi _I\vert \Omega:\Omega \to \Omega _I$ is an
isomorphism, it follows that $\vert \tilde t_I\vert $ is a union
of certain components of $X_I\setminus \Omega _I$. Now, each such
component contains the closed orbit $C_I=G/P_I$.  Therefore
$t_I=0$ and consequently $s=0$, contrary to assumption.
\end {proof}

\subsection {Iwasawa-envelopes}
\label{Iwasawa-envelope}

Here we complete the proof of our main Theorem~\ref{M_D is closed}
by using properties of the Iwasawa-envelope
of the universal domain.  In order to this in $\Omega =G/K$, 
we fix an Iwasawa-Borel subgroup $B$ of $G$ and let $x_0\in \Omega $
be a base point with $\Omega _0:=G_0(x_0)=G_0/K_0$ being the real
symmetric space of basic interest.

Let $H_\Omega $ be the complement in $\Omega $ of the open
$B$-orbit $B(x_0)$, consider the closed $G_0$-invariant set
\begin {equation}\label{Iwasawa complement}
F_\Omega :=\underset{k\in K_0}{\cup}k(H_\Omega )= \underset{g\in
G_0}{\cup}g(H_\Omega) \, ,
\end {equation}
and define the Iwasawa-envelope $\mathcal E_{\mathcal I}(\Omega )$
to be the connected component containing $x_0$ of the complement
of $F_\Omega $ in $\Omega $.  We regard $\mathcal E_\mathcal
I(\Omega )$ as an envelope of $\mathcal U$ in $\Omega $, because
every hypersurface $k(H_\Omega )$ is contained in its complement (\cite
{H}).  In fact, the opposite inclusion also holds (\cite {B}) and
we have the following alternative description of $\mathcal U$
which in fact holds even if $G_0$ is of Hermitian type.
\begin {proposition} \label{Iwasawa-envelope agrees with U}
The Iwasawa-envelope $\mathcal E_\mathcal I(\Omega )$ agrees with
the universal domain $\mathcal U$.
\end {proposition}
Now we give the analogous definition of the Iwasawa-envelope for
an \emph{arbitrary} (algebraic) $G$-compactification $X$ of
$\Omega $. For this let $H_X$ be the closure in $X$ of the
hypersurface $H_\Omega $, and, replacing $H_\Omega $ by $H_X$,
define $F_X$ in the same way as $F_\Omega $.  Then $\mathcal
E_\mathcal I(X)$ is defined to be the connected component
containing $x_0$ of the complement of $F_X$ in $X$.
\begin {theorem} \label {equivalence of Iwasawa}
If $X$ is an arbitrary $G$-equivariant compactification of $\Omega
$, then $\mathcal E_\mathcal I(X)=\mathcal E_\mathcal I(\Omega )$.
\end {theorem}
This follows from Theorem~\ref{no new interior} and 
Proposition~\ref{Iwasawa-envelope agrees with U} in the same way
that Theorem~\ref{M_D is closed} follows from Theorem~\ref{no new interior}
and Theorem~\ref{incidence hypersurfaces}.  But when $G_0$ is not of Hermitian type,
Theorem~\ref{equivalence of Iwasawa} is an immediate consequence of
the following result.
\begin {theorem}\label {Iwasawa in X}
If $X$ is an arbitrary equivariant compactification of $\Omega $
and $G_0$ is not of Hermitian type, then
\begin {equation*}
F_X=\underset{k\in K_0}{\cup}k(H_X)=\underset{g\in
G_o}{\cup}g(H_X)
\end {equation*}
contains the full complement $E=X\setminus \Omega $.
\end {theorem}
Before turning to the proof, let us first
prove a preparatory result which strongly uses the assumption that
$G_0$ is not of Hermitian type.
\begin {lemma} \label{F equals to G/P}
Assume that $G_0$ is not of Hermitian type and let $G/P$ be a
$G$-flag manifold.  If $B$ is an Iwasawa-Borel subgroup of $G$ and
$H$ is the complement of the open $B$-orbit in $G/P$, then
\begin {equation*}
F:=\underset{k\in K_0}{\cup }k(H)=\underset{g\in G_0}{\cup}g(H)
\end {equation*}
is equal to $G/P$.
\end {lemma}
\begin {proof}
Recall that in every open $G_0$--orbit $\gamma $ in $G/P$ there
exists a unique complex $K_0$--orbit $C$. Since $G_0$ is not of
Hermitian type, such \emph{cycles} $C$ are positive-dimensional.

Now the complement of $H$ in $G/P$ is algebraically equivalent to
an affine space $\mathbb C^n$ which contains no positive
dimensional subvarieties.  Thus $H\cap C\not =\emptyset$ for every
base cycle $C$ in every open $G_0$--orbit $\gamma $. The desired
result then follows from the facts that the union of the open
$G_0$--orbits is dense and $F$ is closed.
\end {proof}
The DeConcini-Procesi compactification $X^W$ of $\Omega=G/K$ plays
a special role in the proof of Theorem \ref {Iwasawa in X}. \\

\noindent {\it Proof of Theorem \ref{Iwasawa in X} for $X=X^W$.}
Given $I$ as in $\S$\ref{extending sections}, we show that
$F_{X^W}\supset S_I$. For this let $H$ be the complement of the
open $B$-orbit in $C_I=G/P_I$ as in the above Lemma. By
Proposition~\ref{restriction} (and the brief discussion previous
to it) its pullback $H_I$ to $S_I$ is the zero-set of the
restriction of a $B$-eigensection $t\in \Gamma (X^W,L_I)$
with $\vert t\vert \cap \Omega \not =\emptyset$.  Since $\vert
t\vert \cap \Omega $ is contained in $H_\Omega $, it follows that
$H_I\subset H_{X^W}$.  Now by Lemma \ref{F equals to G/P}
we know that $\underset{k\in K_0}{\cup}k(H)=C_I$.  Thus 
it is immediate that
$\underset{k\in K_0}{\cup}k(H_I)=S_I$. Therefore
\begin {equation*}
F_{X^W}=\underset{k\in K_0}{\cup}k(H_{X^W})\supset \underset{k\in
K_0}{\cup}k(H_I)=S_I\,
\end {equation*}
which completes the proof for $X^W$.\hfill $\square$\\

Now let $X$ be any (algebraic) $G$-equivariant
compactification of $\Omega $. Since it is $G$-equivariantly
birationally equivalent to $X^W$, we have the following diagram:
\begin {gather*}
\begin {matrix}
 & \ \tilde X & \\
\  \ \ \swarrow & &\searrow \\
X^W & & \ \ \ \ X
\end {matrix}
\end {gather*}
Let $\pi :\tilde X\to X^W$ and $p:\tilde X\to X$ denote the
respective proper modifications.\\

\noindent {\it Proof of Theorem \ref{Iwasawa in X} for arbitrary
$X$.}  It is enough to prove this for $\tilde X$, because
$p(F_{\tilde X})=F_X$. Now $\pi :\tilde X\to X^W$ is a
$G$-equivariant proper modification. Since every Borel subgroup
$B$ in $G$ also has an open orbit in $\tilde X$, there are also
only finitely many $G$--orbits in $\tilde X$ and it follows that
the preimage $\tilde S_I=\pi ^{-1}(S_I)$ is also the closure of a
$G$-orbit.

Furthermore, if $H_I=\vert t\vert \cap S_I$ is a $B$-invariant
hypersurface which is defined by a $B$-eigenvector $t\in \Gamma
(X^W,L_I)$ as in the proof for $X^W$, then the
corresponding section $\vert \tilde t\vert \in \Gamma (\tilde
X,\tilde {L}_I)$ of the pullback bundle 
$\tilde {L}_I:=\pi ^*L_I$ has the analogous properties. Namely, the
intersection $\vert \tilde t\vert \cap \tilde \Omega $ of its
support with the open $G$-orbit in $\tilde \Omega $ is nonempty,
and $\vert \tilde t\vert \cap \tilde S_I$ is a $B$-invariant
hypersurface $\tilde H_I$ in $\tilde S_I$.  Thus $\tilde H_I$ is a
subset of the hypersurface $H_{\tilde \Omega }$ which defines the
Iwasawa-envelope in $\tilde \Omega$.

Finally, since the hypersurfaces $k(H_I)$ cover $O_I$ (and
therefore $S_I$) as $k$ runs over $K_0$ and the hypersurfaces
$k(\tilde H_I)$ are just their $\pi $--preimages, it follows that
the hypersurfaces $k(\tilde H_I)$, $k\in K_0$, cover $\tilde S_I$.
This shows that $F_{\tilde X}\supset \tilde S_I$ which completes
the proof in the case of an arbitrary compactification $X$. \hfill
$\square$\\

As a consequence of Theorem~\ref{Iwasawa in X}
we are now able to give the\\

\noindent
{\it Proof of Theorem~\ref{no new interior} in the case where $G_0$ is not
of Hermitian type.} Let $\mathcal V$ be the interior of the closure 
$\cl_X(\mathcal U)$ and observe that $\mathcal V\cap \Omega =\mathcal
U$.  In particular, the intersection $\mathcal V\cap
k(H_X)=\emptyset $ for all $k\in K_0$. On the other hand, if $p\in
\mathcal V\cap E$, it follows from Theorem \ref{Iwasawa in X} that
$p$ is in some such $k(H_X)$ and $\mathcal V \cap k(H_X)\not
=\emptyset$.  Thus no point of $\mathcal V$ is in $E$ and the
first statement of Theorem~\ref{no new interior} follows.  As we explained
directly after the statement of Theorem~\ref{no new interior}, the second statement
is an immediate consequence of the first.  
\hfill \qed\\

As we explained in $\S$\ref{introduction}, Theorem \ref{no new
interior} implies that for every open $G_0$--orbit $D$ of an
arbitrary real form $G_0$ in an arbitrary $G$-flag manifold $Z$
the group theoretically defined cycle space $\mathcal M_D$ is
closed in the full cycle space  
$\mathcal C_q(D)$ and therefore the proof of our main result
is complete.\\

\noindent
{\bf Acknowledgements.} The research for this paper took place
during the 2005 Fukuoka conference on differential geometry and
the authors' visit to Osaka University.  They would like to thank
Professors A.~Fujiki, R.~Kobayashi and Y.~Suyama for making this
possible.  The second author would also like to thank the
Deutsche Forschungsgemeinschaft for longterm support.

Jaehyun Hong\\
Research Institute of Mathematics\\
Seoul National University\\
San 56-1 Shinrim-dong\\
Kwanak-gu\\
Seoul 151-747, Korea\\
jhhong@math.snu.ac.kr\\

\noindent
Alan Huckleberry\\
Fakult\"at und Institut f\"ur Mathematik\\
Ruhr-Universit\"at Bochum\\
D-44780 Bochum, Germany\\
ahuck@cplx.rub.de\\

\end {document}